\documentclass[12pt]{article}
\usepackage[reqno]{amsmath} 
\usepackage{amssymb, amsthm}
\usepackage{enumerate}  
	\usepackage{url} 
	\usepackage{latexsym}
	\usepackage{rotating}
	\usepackage[normalem]{ulem}
	\usepackage{hyperref}
	\usepackage{blkarray}

	\newtheoremstyle{plainsl}%
	{\topsep}
	{\topsep}
	{\slshape} 
	{}
	{\normalfont\bfseries}
	{.}
	{ }
	{}
	
	\swapnumbers
	
	\theoremstyle{plainsl}
	\newtheorem{theorem}{Theorem}[section]
	\newtheorem{lemma}[theorem]{Lemma}

	\newtheorem{proposition}[theorem]{Proposition}
	\normalfont
	\normalfont

	\newcommand\cref[1]{Corollary~\ref{cor:#1}}

	\renewcommand\proof{\noindent\textsl{Proof. }}
	\newcommand\sqr[2]{{\vbox{\hrule height.#2pt
				\hbox{\vrule width.#2pt height#1pt \kern#1pt
					\vrule width.#2pt}\hrule height.#2pt}}}
	\renewcommand\qed{%
		\ifmmode\eqno\sqr53
		\else\nolinebreak\ \hfill\sqr53\medbreak\fi}
	
	\numberwithin{equation}{section}
	
	\newcommand\tr{\text{tr} }

	\makeatletter
	\newcommand\incircbin{\mathpalette\@incircbin}
	\newcommand\@incircbin[2]{\mathbin{\ooalign{\hidewidth$#1#2$\hidewidth\crcr$#1\diamond$}}}
	
	\makeatother

	\newcommand{\ls}[1]
	{\dimen0=\fontdimen6\the\font\lineskip=#1\dimen0
		\advance\lineskip.5\fontdimen5\the\font
		\advance\lineskip-\dimen0
		\lineskiplimit=0.9\lineskip
		\baselineskip=\lineskip
		\advance\baselineskip\dimen0
		\normallineskip\lineskip\normallineskiplimit\lineskiplimit
		\normalbaselineskip\baselineskip
		\ignorespaces}
	
	\usepackage{mathtools}
	
	\usepackage{amsmath}
	
	\newcommand\C{\mathcal{C}}

	\newcommand \qone{\textbf{1}}
	\newcommand \qtwo{\textbf{i}}
	\newcommand \qthree{\textbf{j}}
	\newcommand \qfour{\textbf{k}}
	
	\title{Quantum independence and chromatic numbers}
	\author{Chris Godsil\footnote{Department of Combinatorics and Optimization, University of Waterloo, 200 University Avenue West, Waterloo, ON, Canada N2L 3G1. Email: \protect\url{cgodsil@uwaterloo.ca}. C. Godsil gratefully acknowledges the support of the Natural Sciences and Engineering Council of Canada (NSERC), Grant No. RGPIN-9439.}, 
	Mariia Sobchuk\footnote{Corresponding Author. Department of Combinatorics and Optimization, University of Waterloo, 200 University Avenue West, Waterloo, ON, Canada N2L 3G1. Email: \protect\url{msobchuk@uwaterloo.ca}} }

\begin{document}
		\maketitle
		\begin{abstract}
			We construct a new graph on 120 vertices whose quantum and classical independence numbers are different. At the same time, we construct an infinite family of graphs whose quantum chromatic numbers are smaller than the classical chromatic numbers. Furthermore, we discover the relation to Kochen-Specker sets that characterizes quantum cocliques that are strictly bigger than classical ones. Finally, we prove that for graphs with independence number is two, quantum and classical independence numbers coincide. 
		\end{abstract}
		
	\section{Introduction}
	
	We investigate quantum analogues of the independence number $\alpha_q(X)$ of a graph $X$, and of the chromatic number $\chi(X)$. The quantum independence number was introduced respectively by Man\v{c}inska and Roberson in \cite{RobManQCoc}. The quantum chromatic number appeared first in the work of Avis et al.~\cite{Avis2005AQP}, building on \cite{GalliardWolf},\cite{Cleve2004}.
	
	Define the \textsl{quantum independence number}\label{qcoc} of a graph $X$ to be the largest integer $s$, such that there exists a $|V(X)|\times s$ matrix $P$, such that its entries $P_{vi}$ are projections in $\mathbb{C}^{d\times d}$ and
	\begin{enumerate}[(a)]
		\item $\sum_{v\in V(G)}{P_{vi}}=I_d$. 
		\item ${P_{vi}}{P_{uj}}=0$ for all $i\ne j\in \{1,...,s\}$ and $u\sim_X v$. 
		\item ${P_{vi}}{P_{vj}}=0$ for all $i\ne j\in \{1,...,s\}$ and $v\in V$.
	\end{enumerate}
	We call such matrix $P$ a \textsl{quantum coclique matrix}.
	
	The quantum independence number was also studied in \cite{piovesanthesis}, where, in particular, Teresa Piovesan discovered a graph that has $\alpha<\alpha_q.$ Realising that this graph is a Cayley graph, we find a new graph on 120 vertices with $\alpha<\alpha_q.$ 
	
	We define the \textsl{quantum chromatic number} of a graph $X$ to be the smallest integer $s$, such that there exists a quantum homomorphism from $X$ into $K_s$. This translates into the existence of an $|V(X)|\times t$ matrix $P$, such that its entries $P_{vi}$ are projections in $\mathbb{C}^{d\times d}$ and
	\begin{enumerate}[(a)]
		\item $\sum_{i\in [s]}{P_{vi}}=I_d.$ 
		\item ${P_{vi}}{P_{vj}}=0$ for all $v$ in $V(G)$ and $i\ne j$.
		\item ${P_{vi}}{P_{ui}}=0$ for all $u\sim_X v$.
	\end{enumerate}
	
	The fourteen-vertex graph with $\chi_q<\chi$ is given in \cite{oddities} was recently shown in \cite{lalonde} to be the smallest example where $\chi_q>\chi$. We will show that this graph is an Erd\H{o}s-R\'{e}nyi graph. Next, this observation will lead us to the construction of an infinite family of examples where $\chi>\chi_q$.
	
	Both the quantum independence number and the quantum chromatic number are closely related to the notion of a Projective Kochen-Specker set, defined originally in \cite{scarpathesis}. A \textsl{projective Kochen-Specker set} is a set of $k$ sets, each consisting of $n\times n$ orthogonal projections summing to identity, such that it is impossible to select $k$ pairwise nonorthogonal projections. 
	
	We show that in any graph $\alpha<\alpha_q$ if and only if the entries of quantum coclique matrix form a projective Kochen-Specker set. For chromatic numbers, there already is an analogous result in \cite{scarpathesis}.
	
	In this paper, we include some of the known bounds for $\alpha_q$, which we use to prove that $\alpha=2$ if and only if $\alpha_q=2.$ 
		
	\section{Piovesan's example with $\alpha<\alpha_q$ is Cayley}\label{PiovesansGraph}
	In this section we introduce a graph from \cite{piovesanthesis} whose independence number is strictly smaller than its quantum independence number. We discover that it is also a Cayley graph.
	
	When we say $X$ is an \textsl{orthogonality graph} of a set of vectors $S\in\mathbb{R}^d,$ we mean that $V(X)=S$ and two vectors are adjacent whenever they are orthogonal.
	To begin, the original graph, call it $G_p,$ is the orthogonality graph of the following set of vectors. 
		
		\begin{center}
			\begin{tabular}{ c c c c}
				1: (1,0,0,0) & 2: (0,1,0,0) & 3:(0,0,1,0) & 4: (0,0,0,1)  \\ 
				5: (0,1,1,0) & 6: (1,0,0,-1) & 7: (1,0,0,1) & 8: (0,1,-1,0) \\  
				9: (1,1,1,1) & 10: (1,-1,1,-1) &  11: (1,-1,-1,1)    & 12: (1,1,-1,-1) \\
				13: (1,-1,0,0) & 14: (1,1,0,0) &15: (0,0,1,1) & 16: (0,0,1,-1)\\
				17: (-1,1,1,1) & 18: (1,1,1,-1) & 19: (1,-1,1,1) & 20: (1,1,-1,1) \\
				21: (1,0,1,0) & 22: (0,1,0,1) & 23: (1,0,-1,0) & 24: (0,1,0,-1)
			\end{tabular}
		\end{center}
	This graph was used in the information theoretic context in \cite{leungKS}, and later Piovesan in \cite{piovesanthesis} showed that $\alpha(G_p)<\alpha_q(G_p)$. For convenience, we outline the proof of the latter fact, taken from \cite{piovesanthesis}, below.	
	
	\begin{proposition}
		$\alpha(G_p)<\alpha_q(G_p).$
	\end{proposition}
	\begin{proof}
		Note that each row of the aforementioned table forms a basis of $\mathbb{R}^4$. Therefore, the sets of vertices 
		\begin{align*}
			V_1&=\{(1,0,0,0) ,(0,1,0,0) ,(0,0,1,0) , (0,0,0,1) \},\\
			V_2&=\{ (0,1,1,0) ,(1,0,0,-1) , (1,0,0,1) , (0,1,-1,0)\},\\
			&\vdots \\
			V_6&=\{(1,0,1,0), (0,1,0,1) , (1,0,-1,0), (0,1,0,-1)\}.
		\end{align*}
		partition $V(G)$ into six cliques.
		
		One can check that it is impossible to choose six pairwise nonadjacent elements from each of the 4-cliques, making it also a Kochen-Specker graph and establishing an upper bound on $\alpha$, so $\alpha<6$.
		
		However, it is possible to construct a $|(V(G_p))|\times 6$ quantum coclique matrix $P$ by letting 
		$$P_{ij}=\begin{cases}
			v_iv_i^T \text{ if vertex } v_i \text{ is in clique }j\\
			0\text{ otherwise. }
		\end{cases}$$ Therefore, $\alpha_q\ge 6$.\qed
	\end{proof}
	
	Now we define a Cayley graph. Let $G$ be a group, and $\C$ be an inverse-closed subset of $G$ that does not contain identity. A \textsl{Cayley graph} \cite[Section~3.1]{Godsil2001AlgebraicGT} \textsl{Cay}${(G, \C)}$ is the graph with vertex set consisting of elements of a group $G$ and edge set 
		\[E(\text{Cay}(G,\mathcal{C}))=\{(g,h) : hg^{-1} \in C\}.\]
	Let \[\mathcal{C}:=\{g\in S_4: g^2=()\}.\] 
	
	We provide an explicit isomorphism an isomorphism from $G_p$ into a Cayley graph:
	\begin{theorem}\label{CayPiovesan}
		$G_p$ is isomorphic to $\text{Cay}(S_4,\mathcal{C})$.
	\end{theorem}
	\begin{proof}
			For convenience we denote vertices of $G_p$ not by vectors, but by integers between one and twenty four from the table above. 
		One can check that the following assignment is an isomorphism:
		\begin{center}
			\begin{tabular}{ c c c c}
				1: (1) & 2: (12)(34) & 3:(13)(24) & 4: (14)(23)  \\ 
				5: (23) & 6: (1243) & 7: (1342) & 8: (14) \\  
				9: (124) & 10: (234) &  11: (143)    & 12: (132) \\
				13: (1324) & 14: (1423) &15: (34) & 16: (12)\\
				17: (142) & 18: (134) & 19: (123) & 20: (243)\\
				21: (1432) & 22: (13) & 23: (1234) & 24: (24).\qed
			\end{tabular}
		\end{center}
	\end{proof}

	\section{Construction of the isomorphism to Cayley graph}
	Now we will use quaternions to construct the isomorphism from Theorem \ref{CayPiovesan}. This construction will lead us to a new graph on 120 vertices with $\alpha<\alpha_q$ in the next section.
	
	Before we defined $V_p$ as the unordered set of vertices of $G_p$. Let $V'_p=\{v'_i\}^{24}_i$ be defined as follows:
		\[v'_i=\begin{cases}
			-v_i, \text{ if } i = 7,\,10,\,11,\,13,\,14,\,17,\,18,\,20,\,23\\
			v_i, \text{ otherwise}
		\end{cases}
		\]
		Let $G'_p$ be the orthogonality graph of $V'_p.$
		\begin{lemma}
			$G'_p$ is isomorphic to $G_p.$
		\end{lemma}
		\begin{proof}
			Clearly, for nonzero $\lambda \in \mathbb{R}\setminus\{0\}$, we have $\langle v_i,v_j\rangle = 0$ if and only if $\langle v_i,\lambda v_j\rangle = 0$. Since every $v'_i\in V'_p$  differs from its counterpart $v_i\in V_p$ by at most a multiple of -1, the orthogonality graphs of $V_p$ and $V'_p$ are isomorphic. So $G'_p\cong G_p$. \qed
		\end{proof}
		
		Now we will show how to construct an explicit isomorphism between $G_p$ and $\text{Cay}(S_4,\mathcal{C})$.
		
		\begin{theorem}\label{QuaternionsConstruction}
			$G_p$ is isomorphic to $\text{Cay}(S_4,\mathcal{C})$.
		\end{theorem}
		\begin{proof}
		Let $\mathbb{H}$ denote the quaternions. We denote the basis of $\mathbb{H}$ by $\qone,\qtwo,\qthree,\qfour.$ 
		
		From the previous lemma, $G_p\cong G'_p.$ The idea of the proof is to construct an isomorphism between $\text{Cay}(S_4,\mathcal{C})$ and $G'_p$. First, we will find a relabelling $h: S_4\mapsto \mathbb{H},$ from the vertex set of $G'_p$ into quaternions, distributive over $S_4$ multiplication and mapping the identity permutation to the multiplicative identity of quaternions, $\qone$.
		
		Then we will define a linear map $g$ from $\mathbb{H}\mapsto \mathbb{R}^4$, mapping generators of quaternions to standard basis vectors $\{e_i\}_i$ and distributive over quaternion addition. It will turn out that image under $g$ of $\{h(x)\}_{x\in S_4}$ is just $V(G'_p).$ Moreover, for any $a,b\in S_4,$ it will hold that $\langle g(h(a)),g(h(b))\rangle_{\mathbb{R}^4}=0$ if and only if $ab^{-1}$ is an element of order two in $S_4$. So $g$ is the graph isomorphism between the original Cayley graph, relabelled with $h,$ and $G'_p.$
			
		Let the first four vertices in the embedding in Theorem \ref{CayPiovesan} be denoted by $B_1$, the second four vertices $B_2,$ etc. We observe that the cliques $B_2,...,B_6$ are the cosets of the subgroup $$\{(1),(12)(34),(13)(24),(14)(23)\}$$ of $S_4$, in other words, of the vertices in $B_1$. Denote vertices in $B_1$ by $v_1,...,v_4.$ Choose an assignment as follows. 
			\begin{align*}
				h(v_1)=h(())= q_1,\\
				h(v_2)=h((12)(34)=  q_2,\\
				h(v_3)=h((13)(24))=  q_3, \\
				h(v_4)=h((14)(23))=  q_4.
			\end{align*}
			This is an arbitrary assignment, but we will show that, in fact, $ q_1,...., q_4$ are the basis for quaternions. 
			
			First, we will show that if $i\ne j $, then $q_iq_j \not \in \{ q_1,q_i,q_j\}$. From the above assignment, we have that
			\[q_iq_j=h(v_i)h(v_j)\not \in \{h(v_i),h(v_j)\}=\{ q_1,q_j\}.\]
			Now we show that in order for $h$ to respect $S_4$ multiplication, we must have $h(v_i)h(v_j)\ne  q_1$. Suppose that 
			\begin{align*}
				h(v_i)h(v_j)&=h(v_1), \text{ then multiply on the right by  }h(v_j)\\
				h(v_i)h(v_j)h(v_j)&=h(v_1)h(v_j)=h(v_j) \\
				h(v_i)h(v_j^2)&=h(v_j) \text{, since }h \text{ respects multiplication}\\
				h(v_i)h(())&=h(v_j) \text{, since }v_j^2=1\\
				\pm q_i=q_j,
			\end{align*}
			which is false. Additionally, we will assume that $q_i^2=-q_i$ for $i=2,3,4$.

			The first step towards our goal of mapping elements of $S_4,$ which are vertices of $\text{Cay}(S_4,\mathcal{C}),$ into vectors of $\mathbb{R}^4$ involves representing each vertex as a quaternion. To do so, we will first find images under $h$ of transpositions. Then, we will use multiplication rules on quaternions and the fact that every permutation can be written as a product of transpositions to find images under $h$ of the remaining permutations.  We illustrate the general procedure with one example. The other cases will be analogous.
			\begin{align*}
				h((14))&=a q_1+b q_2+c q_3+d q_4, \quad a,b,c,d\in \mathbb{R},\\
				h((23))&=x q_1+y q_2+zq_3+w q_4, \quad x,y,z,w \in \mathbb{R}.
			\end{align*}

			Using our assumption from the beginning of the proof that $$h((14)(23))= q_4,$$ and realising that the isomorphism $h$ should respect group multiplication, we construct the system of equations:
			\[\begin{cases} 
				h((14))=h((14)(23))h((23)), \\
				h((23))=h((14))h((14)(23)). \\ 
			\end{cases}\]
			
			\[\begin{cases} 
				a q_1+b q_2+cq_3+d q_4= q_4(x q_1+y q_2+zq_3+w q_4), \\
				x q_1+y q_2+zq_3+w q_4=(a q_1+b q_2+cq_3+d q_4) q_4. \\ 
			\end{cases}\]
			Using that 	$ q_4^2=-1:$
			\[\begin{cases} 
				a q_1+b q_2+cq_3+d q_4=x q_4+y q_4 q_2+z q_4q_3-w q_1, \\
				x q_1+y q_2+zq_3+w q_4=a q_4+b q_2 q_4+cq_3 q_4-d q_1. \\ 
			\end{cases}\]
			We obtain:
			\[\begin{cases} 
				a=-w,\\
				d=x ,\\
				cq_3=y q_4 q_2, \\
				b q_2=z q_4q_3,\\
				x=-d ,\\
				w=a ,\\
				y q_2=cq_3 q_4,\\ 
				zq_3=b q_2 q_4.
			\end{cases}\]
			Clearly, $a=d=w=x=0$.
			\[\begin{cases} 
				cq_3=y q_4 q_2, \\
				b q_2=z q_4q_3,\\
				y q_2=cq_3 q_4,\\ 
				zq_3=b q_2 q_4.
			\end{cases}\]
			
			\[\begin{cases} 
				y q_2 q_4=c q_3 q_4 q_4=-cq_3=-y q_4 q_2,\\ 
				zq_3 q_4=b q_2 q_4 q_4=-b q_2=-z q_4q_3,\\
				cq_3 q_2=y q_4 q_2 q_2=-y q_4=-cq_3 q_2 .\\
				
			\end{cases}\]
			Above we established that if $i\ne j $, then $h_ih_j\not \in\{h_1,h_i,h_j\}.$ Thus, we confirm that 
			\begin{align*}
				 q_4q_3&= q_2=-q_3 q_4\\
				 q_2 q_4&=q_3=- q_4 q_2\\
				q_3 q_2&= q_4=- q_2q_3\\
				q_i q_1&= q_1q_i \text { for } i=2,3,4\\
				q_i^2&=- q_1.
			\end{align*}
			We have shown that if $h$ is a group homomorphism from $S_4$ into a group with the property that elements square to minus identity in the image, then the image group is quaternions. In particular, $q_1=\qone, q_2=\qtwo, q_3=\qthree, q_4=\qfour.$
			
			Let $e_1,...,e_4\in \mathbb{R}^4$ denote basis vectors of $\mathbb{R}^4.$
			The signs are arbitrary, and depending on which choice of signs one has, the embedding will be different. Moreover, for convenience we choose integers $\mathbb{Z}_3$ as coefficients, and we have with  
			\begin{align*}
				g(h((14)))=e_2-e_3,\\
				g(h((23)))=e_2+e_3.\\
			\end{align*}
			Similarly, 
			\begin{align*}
				g(h((13)))=e_2+e_4,\\
				g(h((24)))=e_2-e_4,\\
				g(h((12)))=e_3-e_4,\\
				g(h((34)))=e_3+e_4.
			\end{align*}
			To construct representations for all other permutations, write them canonically as a product of transpositions, and use the above-mentioned rules for the multiplication of $q_i$ to get their image under $h$. After applying $g,$ one can check that for any $a,b\in S_4,$ we have $\langle gh(a),gh(b) \rangle =0$ if and only if $ab^{-1}\in\mathcal{C}.$ Thus, we have found an isomorphism from Cay$(S_4,\mathcal{C})$ into $G'_p$.\qed
		\end{proof}
	\section{New graph on 120 vertices with $\alpha <\alpha_q$}
	In Section \ref{PiovesansGraph} we saw that the quantum coclique for $G_p$  was constructed using the fact the $G_p$ was an orthogonality graph, whose vertices conveniently partitioned into cliques. Then, we proved that it was a Cayley graph.
	
	The question arises if the above sequence of observations can be reversed. In other words, can we find more graphs with $\alpha_q>\alpha$ in this way? So we will start with a Cayley graph on $S_n$ with connection set consisting of transpositions in $S_n$, then turn it into an orthogonality graph and hope that $\alpha_q>\alpha$ for such a graph.
	
	In this section we provide one such example of a Cayley graph, refer to it as $G_{120},$ arising from $S_5.$ It holds that $29=\alpha(G_{120})$ and 
	\[\alpha_q(G_{120})\ge 30.\]
	
	We will try to follow the procedure outlined in Theorem \ref{QuaternionsConstruction}. First, we will start with a Cayley graph on $S_5$ with connection set consisting of transpositions in $S_5$. Then we will find a function from $S_5$ to $\mathbb{R}^4.$ Then we use the images of $S_5$ in $\mathbb{R}^4$ to construct an orthogonality graph, and verify that its quantum and classical independence numbers are indeed different.  
	
	To start, assign every vertex of $S_4\subseteq S_5$ the same 4-dimensional vector as in the proof of Theorem \ref{QuaternionsConstruction}. Now we need to decide how the transpositions in $S_5$ but not $S_4$ are represented. Here is the example of this procedure with $(1,5)$. Transpositions $(2,5)$, $(3,5)$ and $(4,5)$ will be analogous.
	
	We want a vector corresponding $(1,5)$ to be adjacent to the vectors corresponding to
	\begin{align*}
		(2,3)=(1,5)\dot (1,5)(2,3),\\
		(2,4)=(1,5)\dot (1,5)(2,4),\\
		(3,4)=(1,5)\dot (1,5)(3,4).
	\end{align*}
	in the orthogonality graph. From the above calculations we had that 
	\begin{align*}
		(2,3)\rightarrow e_2+e_3,\\
		(2,4)\rightarrow e_2-e_4,\\
		(3,4)\rightarrow e_3+e_4.
	\end{align*}
	Now, let $(1,5)$ correspond to $ae_1+be_2+ce_3+de_4$. Then 
	\begin{align*}
		(ae_1+be_2+ce_3+de_4)\cdot (e_2+e_3)=b+c=0,\\
		(ae_1+be_2+ce_3+de_4)\cdot (e_2-e_4)=b-d=0,\\
		(ae_1+be_2+ce_3+de_4)\cdot (e_3+e_4)=c+d=0.
	\end{align*}
	Therefore, $(1,5)$ corresponds to a vector $(*,b,-b,b)$. Let us choose $(0,1,-1,1)$.
	In general, let us use vectors in this representation of $S_5$.
	
	\begin{align*}
		&(15) \text{ to corresponds to a vector }(0,1,-1,1),\\
		&(25) \text{ to corresponds to a vector }(0,1,1,-1),\\
		&(35) \text{ to corresponds to a vector }(0,1,1,1),\\
		&(45) \text{ to corresponds to a vector }(0,-1,1,1).
	\end{align*}
	
	Now, once we have images of all transpositions in $S_5$, we can assign vectors in $\mathbb{R}^4$ to all the remaining elements of $S_5$. We follow a procedure as with $G_p.$ First, write a permutation as a product of transpositions. Then represent each permutation by a vector in $\mathbb{R}^4.$ Then view each vector in $\mathbb{R}^4$ as a sum of quaternions. Multiply two quaternions as usual. In the resulting quaternion, regard the coefficients of $\qone,\mathbf{i},\mathbf{j},\mathbf{k}$ as a vector in $\mathbb{R}^4$. 
	
	Now we will use the vectors assigned to the elements of $S_5$ to construct orthogonality graph $G_{120}$. Turns out, $G_{120}$ can be partitioned into $30$ cliques of size 4, resulting into the lower bound $\alpha_q(G)\ge 30$. Computationally in Sage, one will be able to check that $\alpha(G)=29$. 
	
	The natural question arises: can we apply the method described below to $S_6$ to get a graph with $\alpha<\alpha_q.$ We can not, since we want the new transpositions, for example, $(1,6)$ to be orthogonal to $(2,3)$, $(3,4)$, $(2,4)$, $(1,5)$, $(2,5)$, $(3,5)$, $(4,5).$ Considering which vectors have been assigned to these transpositions earlier, we conclude that $(1,6)$ would have to be a zero vector. So it would be a cone vertex.

	\section{Quantum independence number and Kochen-Specker sets}
	Here we will describe how the quantum coclique matrix determines whether $\alpha<\alpha_q$ or $\alpha=\alpha_q$, without knowing $\alpha$. 
	
	An \textsl{orthogonality graph} of any set of vectors $S$ in an inner product space $\mathcal{H}$ is the graph where the vertices are the vectors in $S$ which are adjacent if the corresponding vectors are orthogonal. We will be working with infinite orthogonality graphs of unit vectors in $\mathbb{C}^d$ and of rank-one projections in $\mathbb{C}^{d}\times \mathbb{C}^d$. The vector inner product will be dot product and the matrix inner product will be the standard trace inner product.  While the two orthogonality graphs are not isomorphic, they are \textsl{homomorphically equivalent}.
	
	\begin{theorem}
		Let orthogonality graph $G$ of the unit vectors in $\mathbb{C}^d$ and $H$ be the orthogonality graph of the rank-one projections in $\mathbb{C}^d\times \mathbb{C}^d$. Then there are homomorphisms $g: V(G)\rightarrow V(H)$ and $h: V(H)\rightarrow V(G)$.
	\end{theorem}
	
	\begin{proof}
		Let  
		\begin{align*} 
			g: V(G)&\rightarrow V(H)\\
			v\mapsto vv^*
		\end{align*}
		For $u,v\in \mathbb{C}^d$ we have that $\langle u, v \rangle =0 $ if and only if $\langle uu^*, vv^* \rangle =0$. It follows that $g$ is a homomorphism.
		Now we claim the following is a homomorphism as well
		\begin{align*} 
			h: V(H)&\rightarrow V(G)\\
			P&\mapsto Pe_1.
		\end{align*}
		Suppose $P_1,P_2$ are the vertices of $H$. Since $\langle P_1,P_2\rangle = \tr P_1P_2 =0$, we have that $P_1P_2=0.$ It follows that $\langle P_1e_1,P_2e_1\rangle=0.$\qed
	\end{proof}
	
	We say a set is a \textsl{transversal for cliques of size $d$} if it contains one vertex from each $d$-clique. A set is \textsl{partial transversal} if it contains at most one vertex from every clique. A function $g: V(G)\rightarrow [0,1]$ is  a \textsl{fractional transversal for $d$-cliques} if $\sum_{x\in C}{g(x)}=1$ for every $d$-clique $C$. 
	
	We denote the orthogonality graph of the unit sphere in $\mathbb{R}^d$ by  $\Omega(d)$. The following result is an immediate consequence of Gleason's theorem:
	\begin{theorem}
		For $d\ge 3$, if $f$ is a fractional transversal on $\Omega(d)$, then there is a matrix $W$ in 
		$\mathbb{C}^d\times \mathbb{C}^d$ such that $f(x)=x^TWx.$\qed
	\end{theorem}
	
	Gleason shows that such $f$ is a continuous map from $\mathbb{R}^d$ into the $[0,1]$ interval. Suppose, we are looking for a binary function that would be a fractional transversal on $\Omega(d)$. Since a binary function only evaluates to $0$ and $1$ it would have to evaluate to 1 on one element of each $d$-clique. Binary functions are discontinuous and so Gleason's theorem proves non-existense of binary fractional transversal. In other words, for the unit sphere there is no binary function that evaluates to one on each clique. 
	
	Now, we can use this result to find chromatic number of the orthogonality graph $\Omega(d)$ of the unit vectors in ${R}^d$. Since there are $d$-cliques, $\chi(\Omega(d))\ge d$. Now, suppose the equality held. It would have to be that each clique has $d$-different colours, which would mean that every colour is present in each clique. Any colour class, say $i$, determines a coclique. Moreover, we could define a characteristic function $f: V(\Omega(d))\rightarrow \{0,1\}$ such that 
	\[f(v)=
	\begin{cases}
		1,&\text{ if } \text{colour}(v)=i; \\
		0,&\text{ otherwise.}
	\end{cases}.\]
	Hence, $f$ evaluates to 1 on each clique, and is therefore a discontinuous fractional transversal, a contradiction. Hence $S_d$ is an infinite graph that is not $d$-colourable. 
	
	\begin{theorem}[De Bruijn-Erd\H{o}s]\cite{ErdosDBrujin}
		An infinite graph $G$ is not $k$-colourable if and only if there exists a finite subgraph 
		of $G$ which is not $k$-colourable.\qed
	\end{theorem}
	
	Therefore, there must be an orthogonality graph $G$ of a finite subset of unit vectors in $\mathbb{R}^d$ that is not $d$-colourable, or, there does not exist a binary function evaluating to one on each clique. Now, if $G$ were not a union of cliques, we would be able to construct $G'$ by adding more vectors to make sure each pair of orthogonal vectors lies in a basis. If $G$ had no strict transversal that was a coclique, $G'$ would not have it as well, because $G$ would be a subgraph of $G'$. We define \textsl{Kochen-Specker graph} to be the orthogonality graph of a union of bases in $\Omega(d).$
	
	However, we will be working mostly with the projective projective Kochen-Specker graphs defined first in \cite{scarpathesis}. A \textsl{projective Kochen-Specker graph} is an orthogonality graph of a set of projections in $\mathbb{C}^d\times \mathbb{C}^d$ such that every transversal contains at least two vertices that are adjacent. The set of vertices of such a graph is called a \textsl{projective Kochen-Specker set}.

\section{Characterization of when $\alpha<\alpha_q$ using Kochen-Specker sets}

	In the previous two sections we saw two graphs that have $\alpha<\alpha_q$. Is there a common pattern to their quantum coclique matrices? This section will answer this question affirmatively. 
	
	Earlier, Scarpa in \cite{scarpathesis} deduced a condition characterising when quantum classical and chromatic numbers differ. While we will cover chromatic numbers later in the paper, we include this result here for reference.
	
	\begin{theorem}\cite{scarpathesis}
		For all graphs $X$, we have that $\chi_q(X)<\chi(X)$ if and only if the orthogonality graph of the entries ${\{P_{vi}\}}_{v\in V(X), i \in [\chi_q]}$ of the $|V(X)|\times c$ quantum coclique matrix $P$ is a projective Kochen-Specker graph.
	\end{theorem}
	
	We prove an analogous statement for the quantum independence number.
	\begin{theorem} 
		For all graphs $X$, we have that $\alpha_q(X)>\alpha(X)$ if and only if the orthogonality graph of the entries ${\{P_{vi}\}}_{v\in V(X), i \in [\alpha_q]}$ of the $|V(X)|\times \alpha_q$ quantum coclique matrix $P$ is a projective Kochen-Specker graph.\qed
	\end{theorem}
	
	\begin{proof}
		For the ease of comprehension, we will sometimes use $k$ to denote $\alpha_q(G)$ in this proof.
		For the forward direction, suppose that $k=\alpha_q(G)>\alpha(G)$. We will prove the contrapositive, that is, if there is a quantum coclique matrix orthogonality graph of whose entries is not a projective Kochen-Specker graph, we will be able to find a classical coclique in $X$ of size $k$.
		
		Suppose $P$ is a $|V(X)|\times k$ quantum coclique matrix with $P_{vi}$ being projections in $\mathbb{C}^d\times \mathbb{C}^d.$ Form an orthogonality graph $X_{\alpha_q}$ with $k|V(X)|$ entries of $P$ as vertices. It is easy to see that $X_{\alpha_q}$ is a union of $k$ cliques, each of size $|V(X)|$, formed by the projections in each column of $P$. By assumption, $X_{\alpha_q}$ is not a projective Kochen-Specker graph. Then there is a transversal $T\subseteq \{P_{vi}\}_{v\in V(X),i\in [k]}$ that is a coclique in $X_{\alpha_q}$.
		
		Now, for each $1\le i\le k$ define $S_i$ to be the set of projections in the $i^{th}$ column of quantum coclique matrix: 
		\[S_i:=\{P_{vi}:v\in V(G)\}.\] 
		Then \begin{align*}
			\sum{S_1}=\sum_{v\in V(G)}{P_{v1}}=I_d \quad &(1),\\
			&\vdots\\
			\sum{S_k}=\sum_{v\in V(G)}{P_{vk}}=I_d \quad &(k).
		\end{align*}
		
		From the definition of quantum coclique each row of the corresponding $|V(G)|\times \alpha_q$ matrix contains pairwise orthogonal projections.
		
		Now, use the transversal $T$ to find a classical $k$-coclique $K$ in the original graph $X$, according to the rule:
		\begin{enumerate}[($\star$)]
			\item $\text {If }(P_{vi}) \in T \text{ for some }i \in [k], \text { include } v \text{ into } K$
		\end{enumerate}
		We will see that $K$ is a classical coclique of size $k$. First, recall from the definition of quantum coclique that if $u$ and $v$ are adjacent, then $P_{vi} P_{uj}=0 \text{ for any }i,j\in[k]$.  Therefore, since $X_{\alpha_q}$ is not a projective Kochen-Specker graph and $T$ is a transversal for $d$-cliques, $T$ contains one projection from each column of $P$, has at most one projection per vertex, and at most one projection will be chosen from any pair of adjacent vertices. Thus, the set of vertices chosen by the $(\star)$ rule, say $K'=\{v_1,...,v_m\}$ is a coclique of size $k$. Now $k\le \alpha(G)\le \alpha_q(G)=k$, so $\alpha(G)=k.$
		
		For the other direction, suppose that $\alpha_q(G)=k$, and the orthogonality graph of entries $\{P_{vi}:v\in V(G),i\in[\alpha_q]\}$  of any quantum coclique matrix is a projective Kochen-Specker graph. Towards a contradiction, we assume that $\alpha(G)=k=\alpha_q(G)$.  Now let $K=\{v_1,...,v_k\}$ be a classical coclique of size $k$ in $G$. Build a $|V(X)|\times k$ quantum coclique matrix as follows. Let the row corresponding to $v_t\in K$ contain $I_d$ in the $t^{th}$ column and zero projections in the remaining entries. For vertices in $G$ that are not in $K$, let the corresponding row consist of zero projections. This is a valid $k$-quantum coclique matrix. Moreover, consider the function 
		$f:S\rightarrow \{0,1\}$
		$$
		f(P_{vi})=
		\begin{cases}
			1, &\text{ if } P_{vi}=I_d;\\
			0, & \text{ otherwise }.\\
		\end{cases}
		$$ 
		Observe that $f$ evaluates to 1 on each column of $P$, and it determines a strict transversal that is a coclique. In the beginning we assumed that the orthogonality graph of $\{P_{vi}\}$ was a projective Kochen-Specker, meaning there can not be a strict transversal that is a coclique. We have arrived at the desired contradiction. \qed
	\end{proof}

	\section{Quantum chromatic number}
	In this section we will define another quantum graph parameter, the quantum chromatic number. In addition, we provide a graph from \cite{oddities} whose quantum chromatic number is strictly smaller than the classical one. In the next section we will generalize this example to an infinite family.
	
	The chromatic number of the graph $X$, denoted $\chi(X)$, is the minimum integer $t$ such that vertices of $X$ can be coloured in $t$ colours with no two vertices of the same colour adjacent. It is well-known that, just like in the case with the independence number, chromatic number can be viewed in terms of homomorphisms.
	\begin{lemma}\cite{Godsil2001AlgebraicGT}
		The \textsl{chromatic number} of a graph $X$ is the least integer $r$, such that there is a homomorphism from $X$ to a complete graph $K_r$. \qed
	\end{lemma}
	
	In the analogy with the classical chromatic number, following \cite{RobersonThesis} we define the quantum independence number of a graph $X$ to be the smallest integer $s$, such that there exists a quantum homomorphism from $X$ into $K_s$. This translates into the following.
	
	For a graph $X$, quantum chromatic number $\chi_q(X)$ is the smallest integer $s$ such that there exists a $|(V(X))|\times s$ matrix $P$, whose entries $P_{ix}$ are projections in $\mathbb{C}^{d\times d}$ such that 
	\begin{itemize}
		\item $\sum_{i\in [s]}{P_{xi}}=I_d.$ 
		\item ${P_{xi}}{P_{xj}}=0\text{ for all }x\in V(G) \text{ and } i\ne j.$
		\item ${P_{xi}}{P_{x'i}}=0\text{ for all }x\sim x'.$ 
	\end{itemize}
	Note that the first condition implies the second. One can also check that if we restrict the definition to the dimension $d=1$, then we just obtain the definition of $\chi$, which suggests:
	
	\[\chi_q\le \chi.\]
	
\begin{figure}[h]
		\centering
		\includegraphics[width=0.6\linewidth]{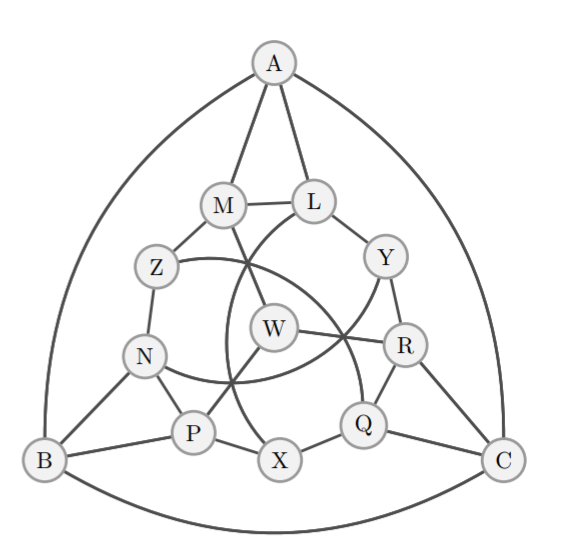}
		\caption{The graph $G_{13}$ from \cite{oddities}}
		\label{fig:G13}
	\end{figure}

	In 2016 Roberson and Man\v{c}inska \cite{oddities} discovered a graph on 14 vertices  for which the quantum chromatic number is less than the classical chromatic number. First, they construct $G_{13}$, a graph defined as an orthogonality graph of columns of the following matrix.
	
	\[
	\left[ \begin{array}{@{}*{13}{c}@{}}
		1 & 0 & 0 & 1 &  1 & 1 & 1  & 0 &  0 & 1 &  1 & -1\\
		0 & 1 & 0 & 1 & -1 & 0 & 0  & 1 &  1 & 1 &  1 &  1\\
		0 & 0 & 1 & 0 &  0 & 1 & -1 & 1 & -1 & 1 & -1 &  1
	\end{array} \right]
	\]
	The illustration of $G_{13}$ is in Figure \ref{fig:G13}.
	
	Then they prove that $G_{14},$ the cone over $G_{13}$, has $\chi(G_{14})>\chi_q(G_{14}).$

	\section{A new infinite family of graphs with $\chi>\chi_q$}
	Now we describe the promised infinite family of graphs with $\chi>\chi_q,$ as a generalisation of the $G_{14}$ example.
	 
	First, an \textsl{Erd\H{o}s-R\'{e}nyi} $ER(p)$ a graph with vertices as basis vectors of one-dimensional subspaces of the vector space over a finite field $\mathbb{F_p}^3$ with two vertices adjacent whenever the corresponding vectors representing subspaces are orthogonal over $\mathbb{F}_p$. 
	
	We realised that $G_{13}$ is an Erd\H{o}s-R\'{e}nyi graph, which lets us generalise Man\v{c}inska-Roberson's construction.  In particular, $G_{13}$ is isomorphic to $ER(3)$. Naturally, vertices of $ER(3)$ can be viewed not only over $\mathbb{Z}_3^3,$ but also over $\mathbb{R}^3$. Also, two vectors representing two vertices of $G_{13}$ are orthogonal when viewed in $\mathbb{R}^3$ if and only if those same vectors viewed over $\mathbb{Z}_3^3,$ are orthogonal. In other words,

	\begin{align*}
		f: V(G_{13})&\rightarrow V(ER(3))\\
		v&\mapsto v
	\end{align*}	
 	is an isomorphism between $V(G_{13})$ and $V(ER(3)).$ But such an embedding map in general may not be an isomorphism between $ER(p)$ and the orthogonality graph of its vertices viewed as vectors in $\mathbb{R}^3$. We must choose representatives for one-dimensional subspaces of $\mathbb{F}_p$ wisely, and then there is a chance. Will show how to define vertices for $ER(p)$ such that the cone of their orthogonality graph over $\mathbb{R}$ will yield a graph where $\chi_q<\chi$.
	
	First, since we are free to choose representatives of subspaces, let us choose the set of vertices of $ER(p)$ to be the $3$-dimensional nonzero vectors over $\mathbb{Z}_p$ with the first nonzero entry being 1:
	\[V':=\{\begin{bmatrix}0,0,1\end{bmatrix},
	\begin{bmatrix}0,1,a\end{bmatrix},
	\begin{bmatrix}1,a,b\end{bmatrix}\}_{a,b\in\{0,...,p-1\}}\setminus\{[0,0,0]\}.\]
	\begin{lemma}\cite{mythesis}
		The vectors in $V'$ belong to distinct one dimensional subspaces of $\mathbb{F}_p^3$. And every one-dimensional subspace of  $\mathbb{F}_p^3$ contains a vector from $V'.$
	\end{lemma}
	
	Now, modify $V'$ as follows. If $x\in V'$ and the $i^{th}$ entry $x_i>\frac{p-1}{2}$ replace $x_i$ with $p-x$. Let the new set be the vertices of $ER(p).$ In this way, the only entries in the vectors of the vertices $V(ER(p))$ will be 
	$$0,\pm 1,...,\pm \frac{p-1}{2}.$$
	
	Now, let any two vertices in $V(ER(p))$ be adjacent whenever they are orthogonal over $\mathbb{R}$ when viewed as vectors over $\mathbb{R}^3.$ Call this graph $ER^\prime(p)$. For example, $ER^\prime(3)$ is just $G_{13}$.
	
	\begin{lemma}
		$ER^\prime(p)$ is connected.
	\end{lemma}
	\begin{proof}
		We will show that between any pair of vertices there exists a path by considering all possible kinds of pairs of vertices. To begin, assume $b\ne 0, d\ne 0$. 
		
		First, consider a pair of vertices of $ER^\prime(p)$ of the form $\begin{bmatrix}0,a,b\end{bmatrix}$ and $\begin{bmatrix}0,c,d\end{bmatrix}$ such that $a,b\in \{0,\pm1,...,\pm \frac{p-1}{2}\}$. Then there is a path consisting of vertices: 
		$$[	0,a,b],\,
		[1,0,0],\,
		[0,c,d].$$
		
		Now look at a different pair of vertices of $ER^\prime(p)$: vertex $v$, labelled by $\begin{bmatrix}0,a,b\end{bmatrix}$ and $u$, labelled by $\begin{bmatrix}1,c,d\end{bmatrix}$ such that $a,b\in \{0,\pm1,...,\pm \frac{p-1}{2}\}$. Then there is a path consisting of vertices: 
		$$[0,a,b],\,
		[1,0,0],\,
		\left[0,1,-\frac{c}{d}\right],\,
		[0,c,d].$$
		
		Finally, between the vertices $[1,a,b]$ and $[1,c,d]$ there is a path 
		$$[1,a,b],\,
		[0,1,-\frac{a}{b}],\,
		[1,0,0],\,
		\left[0,1,-\frac{c}{d}\right],\,
		[1,c,d].$$
		
		Now, in case $b=0,d=0$ there will be a path 
		$$[1,a,0],\,
		[0,0,1],\,
		[1,c,0].$$
		If only one of $b,d$ is $0$, without loss of generality assume $b=0,d\ne 0$, then there will be a path 
		$$[1,a,0],\,
		[0,0,1],\,
		[1,0,0],\,
		\left[0,1,-\frac{c}{d}\right],\,
		[0,c,d].$$
		We have considered all possible pairs of vertices, so we conclude that the graph is connected. \qed
	\end{proof}
	
	\begin{lemma}
		$\chi(ER^\prime(p))\ge 4$.
	\end{lemma}
	\begin{proof}
		Since $ER^\prime(3)$ is the orthogonality graph of the subset of $V(ER^\prime(p))$, $$\chi(ER^\prime(p))\ge \chi(ER^\prime(3))=\chi(G_{13})\ge4.$$
		\qed 
	\end{proof}
	The above lemma shows that the cone over $ER^\prime(p)$ will have $\chi>\chi_q$.
	There are other choices of representatives of subspaces of $ER^\prime(p)$ such that their orthogonality over $\mathbb{R}$ will contain $G_{13}$. For example, as in $G_{13}$ over $\mathbb{Z}_3$, we had vertices $(1,1,0)$ and $(1,-1,0)$ orthogonal over $\mathbb{Z}_3$ and over $\mathbb{R}$. However, if we chose the embedding where the vectors will be not as above, but, say, $(1,1,0)$ and $(1,2,0)$ instead, the orthogonality over $\mathbb{R}$ would have been lost. In general, we can chose a very specific set of entries for the vectors, $0,\pm1,...,\pm \frac{p-1}{2}$, which guarantees that the orthogonality graph over $\mathbb{R}$ will be connected. However, any choice of integers, as long as it includes $0,\pm1$ will work for the proof of the above lemma, but may result in the disconnected graph. Only the connected component containing $G_{13}$ may be considered.
	
	Using similar ideas, we could construct $ER^\prime(p,4)$, with vertices being one-dimensional subspaces of $\mathbb{R}^4$. Technically, it will have no right to be called Erd\H{o}s-R\'{e}nyi any more,  but will contain $G_{14}$ now, the cone over $G_{13}$, and will also yield the a graph with different $\chi$ and $\chi_q$.

	\section{Quantum inertia bound and its application}
	In this section we show the quantum inertia bound technique for bounding the size of quantum coclique. It will be helpful in determining $\alpha_q$ for all graphs with $\alpha=2$.
	
	Let us first recall that the classical \textsl{inertia bound} of $G$. We attribute it to Cvetko\'{v}ic \cite{cvetkovichthesis}. Let  $n_0(G)$, $n_-(G)$, and $n_+(G)$ denote the number of eigenvalues of the adjacency matrix of $G$, which are respectively $0$, strictly negative and strictly positive. Then 
	\[\alpha(G)\le \min\{n_0(G)+n_-(G),n_0(G)+n_+(G)\}.\]
	
	As Wocjan and Elphick have demonstrated, it is also a bound on $\alpha_q:$
	
	\begin{theorem}[Quantum inertia bound]\cite{WocjanElphick}\label{QInertia}
		Let $G$ be a graph, then the following holds:
		\[\alpha_q(G)\le \min\{n_0(G)+n_-(G),n_0(G)+n_+(G)\}.\qed\]
	\end{theorem}
	
	Therefore, if $\alpha(G)$ saturates the inertia bound, then $\alpha(G)=\alpha_q(G):$
	\[\alpha\le \alpha_q(G)\le \min\{n_0(G)+n_-(G),n_0(G)+n_+(G)\}.\qed\]
	
	Now to find the infinite family with $\alpha=\alpha_q,$ we will make use of the following lemma.
	\begin{lemma}\label{SmallerAlpha}
		If a graph $G$ has a quantum-$t$ coclique, it has a quantum-$(t-1)$ coclique.
	\end{lemma}
	
	\proof
	A quantum $(t-1)$-coclique can be obtained from a quantum $t$-coclique matrix by removing any column of projections in the latter.
	
	\begin{lemma}
		If $G$ is a graph on $n\ge2$ vertices, then $\alpha(G)=2$ if and only if $\alpha_q(G)=2$.
	\end{lemma}
	
	\proof
	For the easier only if direction, suppose that $\alpha_q(G)=2$. Since
	\[\alpha(G)\le\alpha_q(G),\] 
	we can only have $\alpha=1$ or $\alpha=2$. Towards a contradiction, suppose it is possible to have $\alpha_q=2$ and $\alpha=1$ in the same graph. The only graph with $\alpha=1$ on $n$ vertices is a complete graph $K_n$. However, $\alpha\le\alpha_q,$ and together with the quantum inertia bound Theorem \ref{QInertia} for $K_n$, we have that 
	\[1=\alpha(K_n)\le \alpha_q(K_n)\le \min\{n_0(K_n)+n_+(K_n),n_0(K_n)+n_-(K_n)\},\]
	For $n\ge 2$, $K_n$ has no 0 eigenvalues and has one positive eigenvalue equal to the degree, so $n_0(K_n)+n_+(K_n)=1$, and thus $\alpha_q(K_n)=1$ as well, which is the contradiction to our initial assumption that $\alpha_q=2$. We conclude that if $\alpha_q(G)=2$, then $\alpha(G)=2$.
	
	Now we will prove that if $G$ is a graph on $n$ vertices with $\alpha(G)=2$, then $\alpha_q(G)=2$. Our strategy is to proceed by contradiction, and assume that it is possible that $\alpha_q(G)\ge 3.$ Note that by Lemma \ref{SmallerAlpha}, if there is a quantum coclique of size $>3$, there is a quantum coclique of size 3. Therefore, without loss of generality, assume that there is a quantum coclique of size 3, i.e. a $|V(X)|\times 3$ matrix $P$ satisfying the quantum coclique definition. Enumerate $V(X)$ as $\{u_1,...,u_n\}$.
	From the quantum coclique definition, for each index (column), we get the column sums:
	\begin{align*}
		P_{u_11}+...+P_{u_n1}&=I,\\
		P_{u_12}+...+P_{u_n2}&=I,\qquad(*)\\
		P_{u_13}+...+P_{u_n3}&=I.
	\end{align*}
	Now, choose an arbitrary vertex $v$, and multiply each of these expressions $	P_{u_1j}+...+P_{u_nj}=I,1\le j\le 3$ by $P_{vi}$ for $i=1,2,3$. When $i=j,$ on the left hand side only the term $P_{vi}$ will survive, so the left hand side will be equal to the right hand side trivially. Hence, we only consider the case when $i\ne j.$ Therefore, the projections $P_{ui}$ with $u\sim v$ are zero. Let $W=\{w_1,...,w_s\},s\le n$ be the non-neighbours of $v$. We are left with the expressions of the following form
	\[P_{vi}(P_{w_1j}+...+P_{w_sj})=P_{vi},\]
	with the conditions that 
	\[v\nsim w_j \forall 1\le j\le s , i\ne j.\]
	
	Moreover, since $\alpha(G)=2$, every set of $3$ vertices of $G$ has at least two vertices with an edge between them. For the set of any two non-neighbours of $v$ and $v$ that looks like $\{w_a,w_b,v\}$, we have to have that there is an edge $w_aw_b$. It follows, that vertices in $W$ forms a clique in $G$.
	Now, multiply each of the following expressions
	\[P_{vi}(P_{w_1j}+...+P_{w_sj})=P_{vi}\]
	for $i\ne j$ and all $i=1,2,3$ on the right by some $P_{w_c,k}$, such that in a such expression  $k\ne i, k\ne j, w_c\in W$. Notice that it is possible to choose such $k$, because the assumption is $\alpha_q\ge 3$. Since $W$ is a clique, any $P_{w_aj}P_{w_ck}=0$. Therefore, the left hand side becomes 0, while on the right hands side we get $P_{vi}P_{w_ck}.$ In this way, we get that $P_{vi}P_{w_ck}=0$ whenever $v$ and $w_c$ are not adjacent and $i\ne k$.  Clearly,  for any $c,w$ we a priori have $P_{vi}P_{wi}=0$ and for adjacent $v\sim w,$ it holds $P_{vi}P_{wi}=0.$ 
	
	Hence, when multiplying all three equations from $(*)$ together, we will be left with 
	\[\sum_{i,j,k=1,...,|V(G)|}{P_{u_i1}P_{u_j2}P_{u_k3}}=0\ne I,\]
	because
	$$P_{u_ia}P_{u_jb}=0\text{ if }u_i\sim u_j$$
	by definition of quantum coclique, and 
	$$P_{u_ia}P_{u_jb}=0\text{ if }u_i\nsim u_j\text{ and }a\ne b$$
	from above. The desired contradiction is achieved.\qed
\section{Acknowledgements}
Both authors acknowledge the support of Chris Godsil’s NSERC (Canada), Grant No. RGPIN-9439.
\bibliographystyle{plain}
\bibliography{MasterPaper2Bib}

	\end{document}